\documentclass[10pt]{amsart}
\usepackage[cp1251]{inputenc}
\usepackage[english,russian]{babel}
\usepackage{amsmath}
\usepackage{amssymb}
\usepackage{amsfonts}

\newtheorem{lemma}{Lemma}
\newtheorem{theorem}{Theorem}

\begin{document}
\renewcommand{\refname}{References}
\renewcommand{\proofname}{Proof.}
\thispagestyle{empty}
\title[Vertex connectivity of Deza graphs]{Vertex connectivity of Deza graphs with parameters of complements to Seidel graphs}
\author{{S.~V.~Goryainov, D.~I.~Panasenko}}%
\address{Sergey Viktorovich Goryainov
\newline\hphantom{iii} Shanghai Jiao Tong University,
\newline\hphantom{iii} 800 Dongchuan RD. Minhang District,
\newline\hphantom{iii} Shanghai, China%
\newline\hphantom{iii} Krasovskii Institute of Mathematics and Mechanics,
\newline\hphantom{iii} S.Kovalevskaya St., 16,
\newline\hphantom{iii} 620990, Yekaterinburg, Russia%
\newline\hphantom{iii} Chelyabinsk State University,
\newline\hphantom{iii} Bratiev Kashirinykh St., 129,
\newline\hphantom{iii} 454001, Chelyabinsk, Russia;}
\email{44g@mail.ru}
\address{Dmitry Igorevich Panasenko
\newline\hphantom{iii} Chelyabinsk State University,
\newline\hphantom{iii} Bratiev Kashirinykh St., 129,
\newline\hphantom{iii} 454001, Chelyabinsk, Russia;}
\email{makare95@mail.ru}

\thanks{\sc Goryainov, S.V., Panasenko, D.I.,
Vertex connectivity of Deza graphs with parameters of complements to Seidel graphs}
\thanks{\copyright \ 2018 Goryainov S.V., Panasenko D.I.,}
\thanks{\rm The work is supported by RFBR according to the projects 16-31-00316 and 17-51-560008}
\thanks{\rm The first author is partially supported by the NSFC (11671258) and STCSM (17690740800)}

\maketitle {\small
\begin{quote}
\noindent{\sc Abstract. }
In this paper we find the vertex connectivity of Deza graphs with parameters of the complements to Seidel graphs.
In particular, we present an infinite family of strictly Deza graphs
whose vertex connectivity is equal to $k-1$, where $k$ is the valency.\medskip

\noindent{\bf Keywords:} strongly regular graph, Deza graph, Seidel graph, vertex connectivity.
 \end{quote}
}

\section{Introduction}
It was proved in \cite{BM85} that the vertex connectivity of a strongly regular is equal to its valency.
In \cite{BK09}, the similar result was obtained in general for distance-regular graphs.
The vertex connectivity of a Cayley graph is at least $2(k+1)/3$, where $k$ is its valency (see \cite{GR01}).
 
In \cite[Theorem 3.1]{EFHHH99}, a construction of strictly Deza graphs was presented.
This construction requires an existence of involutive automorphism of a strongly regular graph $\Gamma$
that interchanges only non-adjacent vertices.
In \cite{GGK14}, the vertex connectivity was studied for strictly Deza graphs obtained from the construction.
It was proved (see \cite[Theorem]{GGK14}) that, if the graph $\Gamma$ has the non-principal eigenvalues $r > 2$ or $s = -2$,
then the vertex connectivity of a Deza graph obtained from the construction is equal to its valency, excepting the case of $(3\times3)$-lattice,
when the valency of the corresponding Deza graph is $4$ and the vertex connectivity is $3$.
The case when the strongly regular graph $\Gamma$ has the non-principal eigenvalue $r \le 2$ was remained open.
In this paper we study the case when $\Gamma$ has the non-principal eigenvalue $r = 1$.

\section{Preliminary}
We consider undirected graphs without loops and multiple edges.

A $k$-regular graph $\Gamma$ on $v$ vertices is called \emph{strongly regular} with parameters $(n,k,\lambda,\mu)$,
if any two vertices $x,y$ in $\Gamma$ have $\lambda$ common neighbours when $x,y$ are adjacent and
$\mu$ common neighbours is $x,y$ are non-adjacent. For a vertex $x$ in a graph $\Gamma$,
the \emph{neighbourhood} $\Gamma(x)$ is the set of all neighbours of $x$ in $\Gamma$.

\begin{lemma}[\cite{BCN89}, Theorem 1.3.1(i)]\label{SRGSpectrum}
Let $\Gamma$ be a strongly regular graph with parameters $(v,k,\lambda,\mu)$, $\mu \ne 0$, $\mu \ne k$.
Then the graph $\Gamma$ has three distinct eigenvalues $k, r, s$, where $k > r > 0 > s$ and the eigenvalues
$r,s$ satisfy the quadratic equation $x^2 + (\mu - \lambda)x + (\mu-k) = 0$.
\end{lemma}

For a graph $\Gamma$, denote by $\overline{\Gamma}$ the complement to $\Gamma$.
\begin{lemma}[\cite{BCN89}, Theorem 1.3.1(x)]\label{SRGComplement}
For a strongly regular graph $\Gamma$  with parameters $(v,k,\lambda,\mu)$,
the complement $\overline{\Gamma}$ is a strongly regular graph with parameters $(n,n-k-1, v-2k+\mu-2, v-2k+\lambda)$
and eigenvalues $n-k-1, -s-1,-r-1$.

\end{lemma}

For any $n \ge 5$, the \emph{triangular graph} $T(n)$
is the graph whose vertices are all $2$-element subsets in $\{1, \ldots, n\}$,
where two vertices are adjacent whenever the intersection of the corresponding subsets has cardinality $1$.
For any $n \ge 3$, the $(n\times n)$-\emph{lattice graph} $L(n)$ is the graph whose vertices are
all ordered pairs of elements from $\{1, \ldots, n\}$, where two vertices are adjacent whenever
the Hamming distance between the corresponding pairs is equal to $1$.
For any $i \in \{1, \ldots, n\}$, the sets $\{(i, j)~|~j \in \{1, \ldots, n\}\}$ and
$\{(j,i)~|~j \in \{1, \ldots, n\}\}$ are called the $i$th \emph{row} and the $i$th \emph{column} of $L(n)$, respectively.

It is easy to show that $T(n)$ is strongly regular with parameters $(n(n-1)/2, 2(n-2), (n-2), 4)$.
and $L(n)$ is strongly regular with parameters $(n^2, 2(n-1), n-2, 2)$.

A strongly regular graph with the eigenvalue $s = -2$ is called a \emph{Seidel graph}.
\begin{lemma}\label{ComplementsToSeidel}
Let $\Gamma$ be a strongly regular graph with the eigenvalue $r = 1$.
Then $\Gamma$ is the complement to a Seidel graph.
\end{lemma}
\proof
It follows immediately from Lemma \ref{SRGComplement}. $\square$

\medskip
A classification of Seidel graphs is given by the following lemma.

\begin{lemma}[\cite{BCN89},Theorem 3.12.4(i)]\label{Seidel}
Let $\Gamma$ be a strongly regular graph with the smallest eigenvalue $s = -2$.
Then $\Gamma$ is a triangular graph $T(n)$ ($n\ge5$), an $(n\times n)$-lattice graph $L(n)$ ($n \ge 3$),
a complete multipartite graph $K_{n\times2}$ ($n\ge 2$) or one of the graphs of Petersen,
Clebsh, Schl{\"a}fli, Shrikhande, or Chang.
\end{lemma}

\medskip
Lemmas \ref{ComplementsToSeidel} and \ref{Seidel} give a classification of strongly regular graphs with the eigenvalue $r = 1$.

A $k$-regular graph $\Delta$ on $v$ vertices is called a \emph{Deza graph} with parameters $(v,k,b,a), b \ge a$,
if any two vertices $x,y$ in $\Delta$ have $a$ or $b$ common neighbours. A Deza graph $\Delta$ is called a \emph{strictly Deza graph},
if the diameter of $\Delta$ is $2$ and $\Delta$ is not strongly regular.
The following lemma gives a construction of strictly Deza graphs.

\begin{lemma}[\cite{EFHHH99},Theorem 3.1]\label{Construction}
Let $\Gamma$ be a strongly regular graph with parameters $(v,k,\lambda,\mu)$, $k\ne\mu$, $\lambda\ne\mu$ and
adjacency matrix $M$. Let $P$ be a permutation matrix that represents an involution $\phi$ of $\Gamma$
that interchanges only non-adjacent vertices. Then $PM$ is the adjacency matrix of a strictly Deza graph $\Delta$
with parameters $(n,k,b,a)$, where $b = max(\lambda,\mu)$ and $a = min(\lambda,\mu)$.
\end{lemma}
An automorphism of a strongly regular graph $\Gamma$ satisfying Lemma \ref{Construction} is called a $\Delta$-\emph{automorphism.}
Since $\phi$ in Lemma \ref{Construction} represents an involution,
the matrix $PM$ is obtained from the matrix $M$ by a permutation of rows in all pairs of rows with indexes $x$ and $y$,
such that $\phi(x) = y$ and $\phi(y) = x$.
Lemma \ref{Images} follows immediately from Lemma \ref{Construction} and shows what is the neighbourhood of a vertex of the graph $\Delta$.
\begin{lemma}\label{Images} For the neighbourhood $\Delta(x)$ of a vertex $x$ of the graph $\Delta$ from Lemma \ref{Construction}, the following conditions hold:
$$
\Delta(x) =
\left\{
  \begin{array}{ll}
    \Gamma(x), & \hbox{if $\phi(x) = x$;} \\
    \Gamma(\phi(x)), & \hbox{if $\phi(x) \ne x$.}
  \end{array}
\right.
$$
\end{lemma}

\medskip

In \cite{GS13}, for the graphs $\overline{L(n)}$ and $\overline{T(n)}$, their $\Delta$-automorphisms were found.

For any $i \in \{1,\ldots,\lfloor n/2\rfloor\}$, let us take the first $i$ pairs of rows in $\overline{L(n)}$
(the pairs of $1$st and $2$nd, $3$rd and $4$th, $\ldots$, $(2i-1)$th and $(2i)$th rows).
Then the permutation that swaps rows in each of the $i$ pairs is a $\Delta$-automorphism of $\overline{L(n)}$. We call such a $\Delta$-automorphism $\overline{L(n)}$ as $i$-\emph{automorphism}.
More precisely, the $i$-automorphism swaps the vertices $(2j,z)$, $(2j+1,z)$ for all $j \in \{1, \ldots, i\}$ and $z \in \{1,\ldots,n\}$.
Note that there are exactly $\lfloor n/2\rfloor$ non-equivalent $i$-automorphisms.
\begin{lemma}[\cite{GS13}, Proposition 6]\label{psiL}
Any $\Delta$-automorphism of $\overline{L(n)}$ is equivalent to an $i$-automorphism for some
$i \in \{1,\ldots,\lfloor n/2\rfloor\}.$
\end{lemma}

\medskip
For any $i \in \{1,2,\ldots,n\}$, denote by $C(i)$ the maximal clique of $T(n)$ induced by
the set of all $2$-subsets that contain $i$. Note that for any distinct $i,j \in \{1,2,\ldots,n\}$
the equality $C(i) \cap C(j) = \{i,j\}$ holds.
The mapping that swaps the vertices $\{1,z\}$ to $\{2,z\}$, for all $z \in {3,4,\ldots,n}$,
is a $\Delta$-automorphism of $\overline{T(n)}$. We call this automorphism as $\{1,2\}$-automorphism.
\begin{lemma}[\cite{GS13}, Proposition 1]\label{psiT}
Any $\Delta$-automorphism of $\overline{T(n)}$ is equivalent to the $\{1,2\}$-automorphism.
\end{lemma}

The \emph{vertex connectivity} $\kappa(\Gamma)$ of a graph $\Gamma$ is the minimum number of vertices whose deletion from
$\Gamma$ disconnects it.
Note that for a $k$-regular graph $\Gamma$ the inequality $\kappa(\Gamma) \le k$ holds.
Let $x$ and $y$ be two vertices of a graph $\Gamma$. Two simple paths connecting $x$ and $y$
are called \emph{disjoint} if they have no common vertices different from $x$ and $y$. A set of vertices $S$
\emph{disconnects} $x$ and $y$ if $x$ and $y$ belong to different connected components of the graph obtained from
$\Gamma$ by deleting $S$. A set $S$
of vertices of a graph $\Gamma$ is called \emph{disconnecting} if it disconnects some two of its vertices. The
following lemma is known as Menger’s theorem.

\begin{lemma}[\cite{H69}, Theorem 5.9]\label{Menger}
The minimum cardinality of a set disconnecting non-adjacent
vertices $x$ and $y$ is equal to the largest number of disjoint paths connecting these vertices.
\end{lemma}

For a path from a vertex $x$ to a vertex $y$, which goes consequently through $m$ vertices $x_1, \ldots, x_m, m \ge 0$,
we use the notation $x \sim x_1 \sim \ldots \sim x_m \sim y$. In the case when $m = 0$ we have the path $x \sim y$ of length $1$,
which is equivalent to the fact that the vertices $x$ and $y$ are adjacent.

In this paper we prove the following two theorems.

In view of Lemma \ref{psiT}, there exists a unique $\Delta$-automorphism of $\overline{T(n)}$,
which is equivalent to the $\{1,2\}$-automorphism. Denote by $\overline{T(n)}'$ the Deza graph
obtained from $T(n)$ with using Lemma \ref{Construction} w.r.t the $\{1,2\}$-automorphism.
\begin{theorem}\label{kappaT}
The vertex connectivity of $\overline{T(n)}'$ is equal to its valency.
\end{theorem}

In view of Lemma \ref{psiL}, any $\Delta$-automorphism of $\overline{L(n)}$ is equivalent to an $i$-automorphism for some
$i \in \{1,\ldots,\lfloor n/2\rfloor\}$. For any $i\in \{1,\ldots,\lfloor n/2\rfloor\}$,
denote by $\overline{L_i(n)}'$ the Deza graph
obtained from $L(n)$ with using Lemma \ref{Construction} w.r.t an $i$-automorphism.
\begin{theorem}\label{kappaL}
For any $i\in \{1,\ldots,\lfloor n/2\rfloor\}$,
the vertex connectivity of $\overline{L_i(n)}'$ is equal to $k-1$, where $k$ is its valency.
\end{theorem}
\section{Proof of Theorem \ref{kappaT}}
By Lemma \ref{SRGComplement}, the graph $\overline{T(n)}$ is strongly regular with parameters $$(n(n - 1)/2, (n^2-5n+6)/2, (n^2-9n+20)/2, (n^2-7n+12)/2).$$
Since, for any $n > 5$, the inequality $(n^2-7n+12)/2 > (n^2-9n+20)/2$ holds, the graph $\overline{T(n)}'$ is a strictly
Deza graph with parameters $$(n(n - 1)/2, (n^2-5n+6)/2, (n^2-7n+12)/2), (n^2-9n+20)/2).$$

Our goal is to prove Theorem \ref{kappaT}. By Lemma \ref{Menger},
it is enough to show that, for any two non-adjacent vertices in $\overline{T(n)}'$, there exists
$(n^2-5n+6)/2$ disjoint paths connecting these vertices.

Through this section, we use the symbols $a,b,c\ldots$ to denote elements of the set $\{1,\ldots, n\}$.
Also, for a $2$-element subset $\{a,b\}$ in $\{1,\ldots, n\}$, we use the shorter notation $ab$.
Two vertices $ab$, $cd$ are adjacent in $\overline{T(n)}$ whenever $|\{a,b\}\cap\{c,d\}| = 0$ holds.
The following lemma shows what is the adjacency rule for the graph $\overline{T(n)}'$.

\begin{lemma}\label{TAdjRule} The following paths of length $1$ occur in $\overline{T(n)}'$.\\
{\rm(1)} $12 \sim ab$, for all $a,b \in \{3,\ldots,n\}, a \not= b$;\\
{\rm(2.1)} $1a \sim 1c$, for all $a,c \in \{3,\ldots,n\}$, $a \not= c$;\\
{\rm(2.2)} $1a \sim bc$, for all $a,b,c \in \{3,\ldots,n\}$, $a \not\in \{b,c\}, b \not= c$;\\
{\rm(3.1)} $2a \sim 2c$, for all $a,b,c \in \{3,\ldots,n\}, a \not= c$;\\
{\rm(3.2)} $2a \sim bc$, for all $a,b,c \in \{3,\ldots,n\}$, $a \not\in \{b,c\}, b \not= c$;\\
{\rm(4)} $ab \sim cd$, for all $a,b,c,d \in \{3,\ldots,n\}$, $\{a,b\} \cap \{c,d\} = \emptyset, a \not= b, c \not= d$.
\end{lemma}
\proof It follows from Lemma \ref{Images} and the definition of the $\{1,2\}$-automorphism. $\square$

\medskip
Now we divide all pairs of non-adjacent vertices in $\overline{T(n)}'$ into several equivalence classes w.r.t. the action of the
$\{1,2\}$-automorphism, and, using Lemma \ref{TAdjRule}, present $(n^2-5n+6)/2$ disjoint paths connecting the vertices in each pair.

(1) Two non-adjacent vertices are both fixed by the $\{1,2\}$-automorphism. Then
these vertices are presented by $ab$ and $bc$, where $a,b,c \in \{3,\ldots,n\}$.
We have the $(n-3)(n-4)/2$ disjoint paths of the form
$$ab \sim de \sim bc, \text{ where }d,e \in \{1,\ldots,n\}, d\ne a,b,c, ~e \ne a,b,c,d;$$
the $n-5$ disjoint paths of the form
$$ab \sim cd \sim be \sim ad \sim bc, \text{ where } d \in \{3,\ldots,n\},~d\ne a,b,c,  
$$$$e \in \{1,\ldots,n\},~e \ne a,b,c,d;$$
the two paths $$ab \sim 1c \sim 1a \sim bc,$$ $$ab \sim 2c \sim 2a \sim bc,$$
which gives a total of $(n^2-5n+6)/2$ disjoint paths.

\medskip
(2) One of two non-adjacent vertices is fixed and another one is moved by the $\{1,2\}$-automorphism.
It is enough to consider the two cases:

(2.1) These two vertices are presented by $ab$ and $1a$, where $a,b \in \{3,\ldots,n\}$, $a\ne b$.
We have the $(n-4)(n-5)/2$ disjoint paths of the form
$$ab \sim cd \sim 1a, \text{ where } c,d \in \{3,\ldots,n\}, c\ne a,b, ~d \ne a,b,c;$$
the $n-4$ disjoint paths of the form
$$ab \sim 1c \sim 1a, \text{ where } c \in \{3,\ldots,n\}, c \ne a,b;$$
the $n-4$ disjoint paths of the form
$$ab \sim 2c \sim 2d \sim bc \sim 1a, \text{ where } c,d \in \{3,\ldots,n\}, c \ne a,b, ~d \ne b,c;$$
the path
$$ab \sim 12 \sim ca \sim 1b \sim 1a, \text{ for some } c \in \{3,\ldots,n\}, c \ne a,b,$$
which gives a total of $(n^2-5n+6)/2$ disjoint paths.

(2.2) These two vertices are presented by $12$ and $1a$, where $a \in \{3,\ldots,n\}$.
We have the $(n-3)(n-4)/2$ disjoint paths of the form
$$12 \sim bc \sim 1a, \text{ where } b,c \in \{3,\ldots,n\}, b\ne a, ~c \ne a,b.$$
Let $\pi$ be a permutation of the set $\{3,\ldots,n\}\setminus\{a\}$ with no fixed points.
Then we have the $n-3$ disjoint paths of the form
$$12 \sim ab \sim 1\pi(b) \sim 1a, \text{ where } b \in \{3,\ldots,n\}, b\ne a.$$
This gives a total of $(n^2-5n+6)/2$ disjoint paths.

For the pairs $ab, 2a$ and $12, 2a$ the arguments are similar.

\medskip
(3) Two non-adjacent vertices are both moved by the $\{1,2\}$-automorphism.
It is enough to consider the following two cases.

(3.1) The two vertices are images of each other w.r.t. the $\{1,2\}$-automorphism,
which means that they are presented by $1a$ and $2a$ for some $a \in \{3,\ldots,n\}$.
We have the $(n-3)(n-4)/2$ disjoint paths of the form
$$1a \sim bc \sim 2a, \text{ where } b,c \in \{3,\ldots,n\}, b\ne a, ~c \ne a,b.$$
Let $\pi$ be a permutation of the set $\{3,\ldots,n\}\setminus\{a\}$ with no fixed points.
Then we have the $n-3$ disjoint paths of the form
$$1a \sim 1b \sim a\pi(b) \sim 2b \sim 2a, \text{ where } b \in \{3,\ldots,n\}, b\ne a.$$
This gives a total of $(n^2-5n+6)/2$ disjoint paths.

(3.2) The two vertices are not images of each other w.r.t. the $\{1,2\}$-automorphism,
which means that they are presented by $1a$ and $2b$ for some $a,b \in \{3,\ldots,n\}, a \ne b$.
We have the $(n-4)(n-5)/2$ disjoint paths of the form
$$1a \sim cd \sim 2b, \text{ where } c,d \in \{3,\ldots,n\}, c\ne a,b, ~d \ne a,b,c.$$
Let us fix some $d\in \{3,\ldots,n\}$ such that $d \ne a,b$.
We have the following three paths:
$$1a \sim 1b \sim ad \sim 2b,$$
$$1a \sim bd \sim 2a \sim 2b,$$
$$1a \sim 1d \sim ab \sim 2d \sim 2b.$$
Let $\pi$ be a permutation of the set $\{3,\ldots,n\}\setminus\{a,b,d\}$ with no fixed points.
Then we have the $n-5$ disjoint paths of the form
$$1a \sim 1c \sim a\pi(c) \sim 2b, \text{ where } c \in \{3,\ldots,n\}, c\ne a,b,d,$$
and the $n-5$ disjoint paths of the form
$$1a \sim bc \sim 2\pi(c) \sim 2b, \text{ where } c \in \{3,\ldots,n\}, c\ne a,b,d.$$
This gives a total of $(n^2-5n+6)/2$ disjoint paths.

The theorem is proved. $\square$

\section{Proof of Theorem \ref{kappaL}}
By Lemma \ref{SRGComplement}, the graph $\overline{L(n)}$ is strongly regular with parameters $$(n^2, (n-1)^2, (n-2)^2, (n-1)(n-2)).$$
Since, for any $n \ge 3$, the inequality $(n-1)(n-2) > (n-2)^2$ holds, the graph $\overline{L_i(n)}'$ is a strictly
Deza graph with parameters $$(n^2, (n-1)^2, (n-1)(n-2), (n-2)^2).$$

Our goal is to prove Theorem \ref{kappaL}. Our approach is the following.
Firstly we prove that, for any $i \in \{1,\ldots,\lfloor n/2\rfloor\}$, the graph $\overline{L_i(n)}'$
has a disconnecting set of size $n^2-2n$.
After that, by Lemma \ref{Menger},
it is enough to show that, for any two non-adjacent vertices in $\overline{L_i(n)}'$, there exists
$n^2-2n$ disjoint paths connecting these vertices.

\begin{lemma} For any $j \in \{1, \ldots, i\}$, the following statements hold.\\
{\rm(1)} The $(2j-1)$th and $(2j)$th rows induce a pair of disjoint cliques in $\overline{L_i(n)}'$.\\
{\rm(2)} The $n^2-2n$ vertices of all rows of $\overline{L_i(n)}'$ but the $(2j-1)$th and $(2j)$th
form a disconnecting set in $\overline{L_i(n)}'$.
\end{lemma}
\proof
(1) It follows from definition of $\overline{L(n)}$, the fact that the $i$-automorphism swaps $(2j-1)$th and $(2j)$th rows,
and Lemma \ref{Images}.

(2) It follows immediately from item (1). $\square$

\medskip
Through this section, we use the symbols $a,b,c\ldots$ to denote elements of the set $\{1,\ldots, n\}$.
Also, for an ordered pair $(a,b)$, where $a,b \in \{1,\ldots, n\}$, we use the shorter notation $ab$.
Two vertices $ab$, $cd$ are adjacent in $\overline{L(n)}$ whenever $a \ne c$ and $b \ne d$ hold.

For any $j \in \{1, \ldots, i\}$, put $(2j-1)' = 2j$ and $(2j)' = 2j-1$.
The following lemma shows what is the adjacency rule for the graph $\overline{L_i(n)}'$.
\begin{lemma}\label{LAdjRule} For any $i \in \{1,\ldots,\lfloor n/2\rfloor\}$,
the following paths of length $1$ occur in $\overline{L_i(n)}'$.\\
{\rm(1)} $ab \sim ac$, for all $a \in \{1, \ldots, 2i\}$, $b,c \in \{1,\ldots,n\}, b \not= c$;\\
{\rm(2)} $ab \sim cd$, for all $a \in \{1, \ldots, 2i\}$, $b, c,d \in \{1,\ldots,n\}, c \not= a', d \ne b$;\\
{\rm(3)} $ab \sim cd$, for all $a \in \{2i+1, \ldots, n\}$, $b, c,d \in \{1,\ldots,n\}, c \not= a, d \ne b$;\\
\end{lemma}
\proof It follows from Lemma \ref{Images} and the definition of the $i$-automorphism. $\square$

\medskip
Now we divide all pairs of non-adjacent vertices in $\overline{L_i(n)}'$ into several equivalence classes w.r.t. the action of the
$i$-automorphism, and, using Lemma \ref{LAdjRule}, present at least $n^2-2n$ disjoint paths connecting the vertices in each pair.

(1) Two non-adjacent vertices are both fixed by the $i$-automorphism.
It is enough to consider the two cases.

(1.1) Two non-adjacent vertices are placed in the same column.
Then
these vertices are presented by $ac$ and $bc$, where $a,b \in \{2i+1, \ldots, n\}$, $a \ne b$,  $c \in \{1,\ldots,n\}$.
We have the $(n-1)(n-2)$ disjoint paths of the form
$$ac \sim de \sim bc, \text{ where }d,e \in \{1,\ldots,n\}, d\ne a,b, ~e \ne c.$$
Let $\pi$ be a permutation of the set $\{1, \ldots, n\} \setminus \{c\}$ with no fixed points.
Then we have the $n-1$ disjoint paths of the form
$$ac \sim bd \sim a\pi(d) \sim bc, \text{ where } d\in \{1,\ldots,n\},  d\ne c.$$
This gives a total of $(n^2-2n+1)$ disjoint paths.

(1.2) Two non-adjacent vertices are placed in the same row.
Then
these vertices are presented by $ab$ and $ac$, where $a\in \{2i+1, \ldots, n\}$,  $b,c \in \{1,\ldots,n\}$, $b \ne c$.
We have the $(n-1)(n-2)$ disjoint paths of the form
$$ab \sim de \sim ac, \text{ where }d,e \in \{1,\ldots,n\}, d\ne a, ~e \ne b,c.$$
Let $\pi$ be a permutation of the set $\{1, \ldots, n\} \setminus \{a\}$ with no fixed points.
Then we have the $n-1$ disjoint paths of the form
$$ab \sim dc \sim \pi(d)b \sim ac, \text{ where } d\in \{1,\ldots,n\},  d\ne a.$$
This gives a total of $(n^2-2n+1)$ disjoint paths.

\medskip
(2) One of two non-adjacent vertices is fixed and another one is moved by the $i$-automorphism.
Then they are placed in the same column, and presented by $ac$ and $bc$, where $a \in \{1,\ldots,2i\}$,
$b\in \{2i+1, \ldots, n\}$, $c \in \{1,\ldots,n\}$.
We have the $(n-1)(n-2)$ disjoint paths of the form
$$ac \sim de \sim bc, \text{ where }d,e \in \{1,\ldots,n\}, d\ne a',b, ~e \ne c.$$
Let $\pi$ be a permutation of the set $\{1, \ldots, n\} \setminus \{c\}$ with no fixed points.
Then we have the $n-1$ disjoint paths of the form
$$ac \sim bd \sim a'\pi(d) \sim bc, \text{ where } d\in \{1,\ldots,n\},  d\ne c.$$
This gives a total of $(n^2-2n+1)$ disjoint paths.

\medskip
(3) Two non-adjacent vertices are both moved by the $i$-automorphism.
It is enough to consider the following cases.

(3.1) Two non-adjacent vertices are placed in rows that are images of each other w.r.t. the $i$-automorphism.

(3.1.1)  Two non-adjacent vertices are placed in the same column.
Then they are presented by $ab$ and $a'b$, where $a,a' \in \{1,\ldots,2i\}$,
$b \in \{1,\ldots,n\}$.
We have the $(n-1)(n-2)$ disjoint paths of the form
$$ab \sim cd \sim a'b, \text{ where }c,d \in \{1,\ldots,n\}, c\ne a,a', ~d \ne b.$$
Take an element $e \in \{1, \ldots, n\} \setminus \{b\}$.
Let $$f:(\{1, \ldots, n\} \setminus \{b,e\}) \rightarrow (\{1, \ldots, n\} \setminus \{a,a'\})$$ be a bijection.
Then we have the $n-2$ disjoint paths of the form
$$ab \sim ac \sim f(c)b \sim a'c \sim a'b, \text{ where } c \in \{1,\ldots,n\}\setminus \{b,e\}.$$
This gives a total of $(n^2-2n)$ disjoint paths.

(3.1.2)  Two non-adjacent vertices are placed in different columns.
Then they are presented by $ab$ and $a'c$, where $a,a' \in \{1,\ldots,2i\}$,
$b,c \in \{1,\ldots,n\}$, $b \ne c$.
We have the $(n-2)^2$ disjoint paths of the form
$$ab \sim de \sim a'c, \text{ where }d,e \in \{1,\ldots,n\}, d\ne a,a', ~e \ne b,c.$$
Let $$f:(\{1, \ldots, n\} \setminus \{b,c\}) \rightarrow (\{1, \ldots, n\} \setminus \{a,a'\}),$$
$$g: (\{1, \ldots, n\} \setminus \{a,a'\}) \rightarrow (\{1, \ldots, n\} \setminus \{b,c\}),$$
 be two bijections.
Then we have the $n-2$ disjoint paths of the form
$$ab \sim ad \sim f(d)c \sim a'c, \text{ where } d \in \{1,\ldots,n\}\setminus \{b,c\}$$
and the $n-2$ disjoint paths of the form
$$ab \sim ec \sim a'g(e) \sim a'c, \text{ where } e \in \{1,\ldots,n\}\setminus \{a,a'\}.$$
This gives a total of $(n^2-2n)$ disjoint paths.

(3.2) Two non-adjacent vertices are placed in rows that are not images of each other w.r.t. the $i$-automorphism.
Then they are placed in the same column, and presented by $ac$ and $bc$, where $a,b \in \{1,\ldots,2i\}, a'\ne b, c \in \{1,\ldots,n\}$.
We have the $(n-1)(n-2)$ disjoint paths of the form
$$ac \sim de \sim bc, \text{ where }d \in \{1,\ldots,n\}, d\ne a',b', ~e \ne c.$$
Let $\pi$ be a permutation of the set $\{1, \ldots, n\} \setminus \{c\}$ with no fixed points.
Then we have the $n-1$ disjoint paths of the form
$$ac \sim b'd \sim a'\pi(d) \sim bc, \text{ where } d \in \{1,\ldots,n\}\setminus \{c\}.$$
This gives a total of $(n^2-2n+1)$ disjoint paths.

The theorem is proved. $\square$

\section{Concluding remarks}
We have found the vertex connectivity of the strictly Deza graphs obtained from 
the complements to triangular and $(n\times n)$-lattice graphs. Let us make some
remarks on the other graphs given by Lemma \ref{Seidel}. Since the complement to 
a complete multipartite graph $K_{n\times2}$ is a disjoint union of edges, we have nothing to prove.
An exhaustive computer search shows that the complements to Shrikhande graph and one of three Chang graphs have 
no $\Delta$-automorphisms; the complements to Schl{\"a}fli graphs and two of three Chang graphs
have a unique $\Delta$-automorphism; the complement to Clebsh graph has precisely two non-equivalent $\Delta$-automorphisms.
For all strictly Deza graphs obtained from Lemma \ref{Construction} w.r.t. these $\Delta$-automorphisms, 
their vertex connectivity equals to the valency.

Strongly regular graphs with $r = 2$ were studied in \cite{KMP10}. 
There were presented three infinite families of possible parameter tuples for such graphs.
For infinitely many of those parameter tuples an existence of a graph with such parameters is unknown.

Note that, for any even $n$, $n \ge 4$, the graph $\overline{L_{\frac{n}{2}}(n)}'$ can be regarded as a Cayley graph
of the group $\mathbb{Z}_n \times \mathbb{Z}_n$, which gives an infinite family of Cayley-Deza graphs
whose vertex connectivity is equal to $k-1$, where $k$ is the valency. 

In \cite{GIKMS17}, strictly Deza graphs with disconnected second neighbourhood of a vertex
were studied. It was proved that if all second neighbourhoods of vertices in a strictly Deza graph $\Delta$ are disconnected,
then $\Delta$ is either edge-regular or co-edge regular.
Let us notice that the second neighbourhood of the vertex $\{1,2\}$ in the graph $\overline{T(n)}'$,
which is neither edge-regular nor co-edge-regular,
is a disjoint union of two cliques induced by the sets $\{\{1,j\}~|~j\in\{3,\ldots,n\}\}$ and $\{\{2,j\}~|~j\in\{3,\ldots,n\}\}$.
Excepting the vertex $\{1,2\}$, the second neighbourhoods of all vertices in $\overline{T(n)}'$ are connected.

\end{document}